\documentclass[12pt,reqno]{amsart}
\usepackage{hyperref}
\usepackage{amssymb,amsmath,amsfonts,latexsym}
\usepackage{bm}
\usepackage{mathtools,mathrsfs}
\usepackage{enumerate}
\usepackage{array,graphics,color}
\allowdisplaybreaks

\numberwithin{equation}{section}
\newtheorem{dfn}{Definition}[section]
\newtheorem{thm}[dfn]{Theorem}
\newtheorem{que}[dfn]{Question}
\newtheorem{lma}[dfn]{Lemma}

\newtheorem{ppsn}[dfn]{Proposition}
\newtheorem{crlre}[dfn]{Corollary}

\newtheorem{rmrk}[dfn]{Remark}

\DeclarePairedDelimiterX{\norm}[1]{\lVert}{\rVert}{#1}
\DeclarePairedDelimiterX{\bnorm}[1]{\big\lVert}{\big\rVert}{#1}
\DeclarePairedDelimiterX{\Bnorm}[1]{\Big\lVert}{\Big\rVert}{#1}
\newcommand\at[2]{\left.#1\right|_{#2}}

\newcommand{\R}{\mathbb{R}}

\newcommand{\Nat}{\mathbb{N}}

\newcommand{\blm}{\mathcal{B}(\ell_2^m)}

\newcommand{\Tr}{\operatorname{Tr}}
\newcommand{\sgn}{\operatorname{sgn}}

\newcommand{\ddt}{\dfrac{d}{dt}}

\begin{document}

	\title[On Isometric Embeddability of non-commutative Quasi-Banach spaces]{On Isometric Embeddability of $S_q^m$ into $S_p^n$ as non-commutative Quasi-Banach spaces}
	
	\author[Chattopadhyay] {Arup Chattopadhyay}
	\address{Department of Mathematics, Indian Institute of Technology Guwahati, Guwahati, 781039, India}
	\email{arupchatt@iitg.ac.in, 2003arupchattopadhyay@gmail.com}
	
	\author[Hong] {Guixiang Hong}
	\address{Institute for Advanced Study in Mathematics, Harbin Institute of Technology, Harbin 150001, China}
	\email{ {\color{red}gxhong@hit.edu.cn}}

	\author[Pradhan] {Chandan Pradhan}
	\address{Department of Mathematics, Indian Institute of Technology Guwahati, Guwahati, 781039, India}
	\email{chandan.math@iitg.ac.in, chandan.pradhan2108@gmail.com}
	
	\author[Ray] {Samya Kumar Ray}
	\address{Stat-Math Unit, Indian Statistical Institute, 203, B. T. Road, Kolkata 700108, India}
	\email{samyaray7777@gmail.com}

	

	\begin{abstract}
		The existence of isometric embedding of $S_q^m$ into $S_p^n$, where $1\leq p\neq q\leq \infty$ and $m,n\geq 2$ has been recently studied in \cite{JFA22}. In this article, we extend the study of isometric embeddability beyond the above mentioned range of $p$ and $q$. More precisely, we show that there is no isometric embedding of the commutative quasi-Banach space $\ell_q^m(\R)$ into $\ell_p^n(\R)$, where $(q,p)\in (0,\infty)\times (0,1)$ and $p\neq q$. As non-commutative quasi-Banach spaces, we show that there is no isometric embedding of $S_q^m$ into $S_p^n$, where $(q,p)\in (0,2)\setminus \{1\}\times (0,1)$ $\cup\, \{1\}\times (0,1)\setminus \{\frac{1}{n}:n\in\mathbb{N}\}$ $\cup\, \{\infty\}\times (0,1)\setminus \{\frac{1}{n}:n\in\mathbb{N}\}$ and $p\neq q$. Moreover, in some restrictive cases, we also show that there is no isometric embedding of $S_q^m$ into $S_p^n$, where $(q,p)\in [2, \infty)\times (0,1)$. A new tool in our paper is the non-commutative Clarkson's inequality for Schatten class operators. Other tools involved are the Kato-Rellich theorem and multiple operator integrals in perturbation theory, followed by intricate computations involving power-series analysis.
	\end{abstract}
	\maketitle
	
	\section{Introduction and main results}
	In a recent article \cite{JFA22}, the authors have studied isometric embeddability of $S_q^m$ into $S_p^n$ for $1\leq p\neq q\leq\infty$ and $m,n\geq 2.$ This partially established non-commutative analogues of the results appearing in \cite{Ly08}, \cite{Ly09}, \cite{LySh01}, \cite{LySh04}, \cite{LyVa93}. The motivation behind this study was two-fold. First of all the study of isometric embeddability for commutative $L_p$-spaces is an exciting topic of research. This began from the seminal paper of Banach \cite{Ba32} and continued by several authors. We refer to \cite{DeJaPe98}, \cite{FlJa03} and \cite{FlJa07} and references therein for a comprehensive study. We also point out \cite{He71}, \cite{KoKo01}, \cite{Le37}, \cite{Ro73}, \cite{Sc38} for exciting connections with various other fields of mathematics for example probability theory and combinatorics. Motivated by quantum mechanics, non-commutative analysis has made great advancements in the past few decades. The theory of non-commutative $L_p$-spaces is an outgrowth of this direction of research. In \cite{Ar75}, \cite{Ye81}, \cite{JuRuSh05}, the authors have studied isometries on non-commutative $L_p$-spaces and established non-commutative analogues of results of Banach \cite{Ba32} and Lamperti \cite{La58}. A remarkable progress took place after the successful development of operator space theory and non-commutative probability theory. We refer to \cite{JuPa082}, \cite{Pi04}, \cite{Ju05}, \cite{Xu06}, \cite{Ju00}. We also refer to \cite{JuPa10}, \cite{JuPa08}, \cite{JuPa082}, \cite{JuNiRuXu04}, \cite{SuXu03}, \cite{JuRuSh05}, \cite{RaXu03}, \cite{HaaRoSu}, \cite{JuSuZa17} and references therein for more information in this direction of research. Therefore, in view of the results appearing in \cite{Ly08}, \cite{Ly09}, \cite{LySh01}, \cite{LySh04}, \cite{LyVa93}, it is natural to study non-commutative analogues of them in the context of non-commutative isometric or isomorphic embedding theory. The second motivation behind \cite{JFA22} was again operator space theory and its connections to boundary normal dilation, where such non-commutative isometric embeddability was crucial. This was also motivated from a guess of Pisier and Misra. We refer to \cite{Ra20} for more on this (see also \cite{GuRe18}).
	
	Beyond the range of $(p,q)\in[1,\infty]\times[1,\infty],$ some results are also known. By \cite{Ju00} for all $0<p<q<2,$ $L_q(\mathcal R\overline{\otimes}B(\ell_2))$ and hence $L_q(\mathcal R)$ embeds isometrically into $L_p(\mathcal R)$ where $\mathcal R$ is a hyperfinite type $\text{II}_1$ factor. This generalized remarkable results due to \cite{BrDiKr66} in the non-commutative setting. In \cite{SuXu03}, Sukochev and Xu studied when $L_p(\mathcal N)$ embeds into $L_p(\mathcal M)$ for $0<p<1,$ where $\mathcal M$ and $\mathcal N$ are semifinite von Neumann algebras. We also refer \cite{RaXu03} and \cite{Ra08} for related work. Despite these remarkable developments, the case when $S_q^m$ embeds isometrically into $S_p^n$ is not well studied for $0<p\neq q\leq\infty$. The question when $S_q^m$ embeds isometrically into $S_p^n$ for $0<p\neq q\leq\infty$, was asked by Q. Xu to the authors. 
	
	\begin{que}\label{Xu}
		Let $0<p\neq q\leq\infty$ and $m,n\geq 2.$ When does $S_q^m$ embeds isometrically into $S_p^n$?
	\end{que} 
	
	Note that \cite{JFA22} and \cite{Ra20} give partial answers to the Question \ref{Xu} when $1\leq p\neq q\leq \infty.$ The main novelty of our work in this article is to study isometric embeddability of $S_q^m$ into $S_p^n$ where $0<p\neq q\leq\infty$ as non-commutative quasi-Banach spaces which goes beyond the foregoing range of $p$ and $q$. We also study the classical case, i.e.  when $\ell_q^m$ embeds isometrically into $\ell_p^n$ for $m,n\geq 2$ and $0<p\neq q\leq\infty.$ Note that the case when $1\leq p\neq q\leq\infty$ have been extensively studied in \cite{Ly08}, \cite{Ly09}, \cite{LySh01}, \cite{LySh04}, \cite{LyVa93}. However, the authors did not consider the cases when $p$ and $q$ are allowed to be strictly less than $1$ and to our surprise there is no more study in the literature on this. Thus our study is new even in the commutative case. Our first theorem is the following. 
	\begin{thm}[Isometric Embeddability between commutative quasi-Banach spaces $\ell_p^n(\mathbb{K})$]\label{main1}
		Let $(q,p)\in (0,\infty)\times (0,1)$, and $2\leq m\leq n<\infty$ and $p\neq q$. Then there is no isometric embedding of $\ell_q^m(\mathbb{K})$ into $\ell_p^n(\mathbb{K})$ for each of the following cases:
		\begin{enumerate}
			\item $\mathbb{K}=\mathbb{C}$ and $q\in (0,\infty)\setminus 2\mathbb{N}$.
			\item $\mathbb{K}=\R$.
		\end{enumerate} 
	\end{thm}
	In the above theorem, the results for real and complex cases are \textsf{different}. The technical reason behind this is that if we consider the complex field instead of $\R$, that is $\ell_q^m(\mathbb{C})$, then in the proof of Theorem \ref{main1}, it may happen that $\at{\frac{d^2}{dt^2}}{t=0}\norm{\bm{a}+t\bm{b}}_p^p=0$, where $\|\bm{a}\|_p=\|\bm{b}\|_p=1$ and in that case we can't conclude anything. Please see Remark \ref{ramk3.1} for more on this. These kind of facts are one of the many subtleties in our work.
	\begin{thm}[Isometric Embeddability between non-commutative quasi-Banach spaces $S_p^n$]\label{main2}
		Let $p\in (0,1)$, and $2\leq m\leq n$. 
		\begin{enumerate}
			\item There is no isometric embedding of $S_q^m$ into $S_p^n$ for $q\in(0,2)\setminus\{1\}$ and $p\neq q$. 
			\item Let $2\leq q< \infty.$ Then there is no isometric embedding $T:S_q^m\to S_p^n$ with\\ $T(\text{diag}(1,0,\ldots,0))=A,$ and $T(\text{diag}(0,1,\ldots,0))=B$ such that 
			\begin{itemize}
				\item  $A,B\in M_n^{sa},$ 
				\item either $A\geq 0$ or $A\leq 0$.
			\end{itemize}
			\item There is no isometric embedding of $S_1^m$ into $S_p^n$ for $p\in(0,1)\setminus\{\frac{1}{k}:k\in\mathbb{N}\}$. 
			\item There is no isometric embedding of $S_\infty^m$ into $S_p^n$ for $p\in(0,1)\setminus\{\frac{1}{k}:k\in\mathbb{N}\}$. 
		\end{enumerate}
	\end{thm}
	We refer to Section \ref{pre} for unexplained notations in the above theorems. Note that Theorem \ref{main2} together with results from \cite{JFA22} and \cite{Ra20} give an answer to the Question \ref{Xu} for a wide range of $p$ and $q.$ \textsf{It is interesting to observe that in view of \cite{Ju00}, Theorem \ref{main2} indicates sharp contrast between finite dimensional and infinite dimensional cases and shows that isometric embeddability problem between non-commutative $L_p$-spaces also depends on the type of the underlying von Neumann algebras.} Each part of Theorem \ref{main2} has been proven using different methods and each of them requires specific techniques developed before that. For these reasons, we have proved them separately in Theorem \ref{th1}, Theorem \ref{th2}, Theorem \ref{th3}, Theorem \ref{th4} and Theorem \ref{th5}.

	The paper is organized as the following. In Section \ref{pre}, we recall necessary background and  prove some useful results. In Section \ref{comiso}, we study the isometric embeddability of $\ell_q^m$ into $\ell_p^n$ as commutative quasi-Banach spaces. In Section \ref{ncomiso}, we recall and prove many facts about operator derivatives and study the isometric embeddability of $S_q^m$ into $S_p^n$ as non-commutative quasi-Banach spaces.
	\section{Preliminaries:}\label{pre}
	We let $\mathbb{K}$ denote the scalar field. That is $\mathbb{K}=\R $ or $\mathbb{C}$. For $0<p<\infty$, we denote by $\ell_p^n(\mathbb{K})$ the $\mathbb{K}$-linear space $\mathbb{K}^n$ equipped with the $\ell_p$ semi-norm
	
	\[\|\bm{a}\|_p=\left(\,\sum_{k=1}^{n}|a_k|^p\,\right)^{\frac{1}{p}}.\]
	
 For simplicity, we write $\ell_p^n(\mathbb{C})=\ell_p^n$. Let $\mathcal M$ be a von Neumann algebra with a normal semifinite faithful trace $\tau.$ Let $\mathcal{S}(\mathcal{M})$ be the linear span of positive elements with finite support. Let $0<p<\infty.$ For $x\in\mathcal{S}(\mathcal{M}),$ define $\|x\|_p:=(\tau(|x|^p))^{\frac{1}{p}}.$ For $p\geq 1,$ $\mathcal{S}(\mathcal{M})$ is a normed space and for $0<p<1,$ it is a quasi-normed space. Then $L_p(\mathcal M)$ is defined to be $\overline{(\mathcal{S}(\mathcal M),\|.\|_p)}$ with respect to the metric $d_p(x,y):=\|x-y\|_p^p.$ One denotes $L_\infty(\mathcal M)=\mathcal M.$ Thus for $1\leq p\leq \infty,$ $L_p(\mathcal M)$ is a Banach space and for $0<p<1,$ it becomes a quasi-Banach space. When $\mathcal M=\mathcal B(\ell_2^n)$ with the usual trace $Tr,$ the corresponding non-commutative $L_p$-space is called  the Schatten $p$-class and denoted by $S_p^n$. The set of all $n\times n$ complex matrices is denoted by $M_n$. We denote $M_n^{sa}$ to be the set of all $n\times n$ self-adjoint matrices. For $0<p<1$, it is well known that $\|.\|_{L_p(\mathcal M)}$ satisfies
	\begin{enumerate}
		\label{sn}\item $\|A+B\|_p\leq 2^{\frac{1}{p}-1} \left(\|A\|_p+\|B\|_p \right),\quad\forall\, A,\,B\in L_p(\mathcal M),$
		\item $\|UAV\|_p=\|A\|_p$ for all unitary operators $U, V$ in $\mathcal M$, and $A\in L_p(\mathcal M)$, that is $\|\cdot\|_p$ is a unitary invariant quasi-norm.
	\end{enumerate}
	It is natural to ask that for which class of operators in $S_p^n$ satisfies the triangular inequality (i.e. $\|A+B\|_p\leq \|A\|_p+\|B\|_p $). The problem of the isometric embeddability of $S_q^m$ into $S_p^n$ for $1\leq q\leq \infty$ can possibly shed some light on this question.
 We now recall a few useful tools.
	\subsection{ Non-commutative Clarkson's inequality for Schatten class operators:} A new tool in our work is the non-commutative Clarkson's inequality. This helps us to obtain the crucial Lemma \ref{lem2.2}. However, we want to emphasize that non-commutative Clarkson's inequality alone cannot solve our problem as we know that non-commutative Clarkson's inequality is also true for general non-commutative $L_p$-spaces but the existence of isometric embeddability between non-commutative $L_p$-spaces are completely different if the underlying von Neumann algebras are of different type as we have pointed out before.
	We refer \cite{Mc67} for the following inequality.
	\begin{thm}
		Let $A,B\in S_p^n$. Then for $0\leq p\leq 2$, $A, B$ satisfies 
		\begin{align}\label{cl}
			 \|A-B\|_p^p + \|A+B\|_p^p\leq 2\left(\|A\|_p^p+\|B\|_p^p\right),
		\end{align}
	and for $2\leq p<\infty$, 
	\begin{align}\label{cll}
		 \|A-B\|_p^p + \|A+B\|_p^p\geq 2\left(\|A\|_p^p+\|B\|_p^p\right).
	\end{align}
	\end{thm}

		Observe  that, for $0<q<\infty$, $0<p< \infty$, $2\leq m\leq n<\infty$, isometric embeddability of $S_q^m$ into $S_p^n$ implies existence of a linear isometry  $T:\ell_q^2(\mathbb{C})\to S_p^n $, which further implies the following equality
		\begin{align}\label{equality}
			\|(1,z)\|_q=\|A+zB\|_p \text{ for all } z\in\mathbb{C}, \text{ where } T((1,0)):=A, \text{ and } T((0,1)):=B.
		\end{align}
		Now if we assume that for $0<q<\infty,~ 0< p\leq 2 $, $S_q^m$ isometrically embeds into $S_p^n$, then non-commutative Clarkson's inequality \eqref{cl} and the above identity \eqref{equality} together implies
		\begin{align}\label{cleq}
			\|(1,z)\|_q\leq \|(1,z)\|_p \text{ for all } z\in\mathbb{C}.
		\end{align} In \eqref{cleq}, the last inequality implies that $p\leq q$. On the other hand if we assume that for $0<q<\infty, 2<p< \infty$, $S_q^m$ isometrically embeds into $S_p^n$, then non-commutative Clarkson's inequality \eqref{cl} and the above identity \eqref{equality} together implies
	\begin{align}\label{cleqq}
	\|(1,z)\|_q\geq \|(1,z)\|_p \text{ for all } z\in\mathbb{C}.
\end{align}
	  Therefore, in this case, $p\geq q$ is necessary for isometric embeddability of $S_q^m$ into $S_p^n$. But, in both cases, non-commutative Clarkson's inequality does not give anything more to these embedding problems. In the following lemma, we summarize the above discussions.
	  \begin{lma}\label{lem2.2}
	  	Let $2\leq m\leq n<\infty,\, 0<p\neq q<\infty$. Suppose $S_q^m$ is isometrically embeds into $S_p^n$. 
	  	\begin{enumerate}
	  		\item If $0<q<\infty$, $0<p\leq 2$, then $p<q$.
	  		\item If $0<q<\infty$, $2\leq p<\infty$, then $p>q$.
	  	\end{enumerate}
	  \end{lma}
\begin{rmrk}
Let $\mathcal M$ be the hyperfinite type $II_1$ factor. Note that as in Lemma \ref{lem2.2}, one can easily see that if $0<q<\infty$ and $0<p\leq 2$, then $L_q(\mathcal M)$ cannot embed isometrically into $L_p(\mathcal M)$ whenever $p>q.$ Also for $0<q<\infty$ and $2\leq p<\infty,$ $L_q(\mathcal M)$ cannot embed isometrically into $L_p(\mathcal M)$ if $p<q.$
\end{rmrk}

	   To continue the study of the isometric embedding problem, here we need to introduce two more novel ingredients, which were also crucial in our previous study \cite{JFA22}, namely the Kato-Rellich theorem and multiple operator integrals in perturbation theory.

	\subsection{A simple version of Kato-Rellich theorem:}
	\begin{thm}$($see \cite[Page 31, Chapter I]{ReBe69}, \cite[P. 21, Theorem 1]{Ba85}$)$\label{kato} Let $A,B\in M_n^{sa}.$ Then, for all $t_0\in\mathbb R$, there exists $\epsilon(t_0)>0$ and real-analytic functions $\lambda_k:(t_0-\epsilon(t_0),t_0+\epsilon(t_0)) \to \mathbb R,$ where $1\leq k\leq n$ and real-analytic functions $u_{ij}:(t_0-\epsilon(t_0),t_0+\epsilon(t_0)) \to \mathbb R,$ where $1\leq i,j\leq n$ such that
		\begin{itemize}
			\item[(i)]for all $t\in(t_0-\epsilon(t_0),t_0+\epsilon(t_0))$, $\{\lambda_1(t),\dots,\lambda_n(t)\}$ is the complete set of eigenvalues of $A+tB$ counting multiplicity.
			\item[(ii)] for all  $t\in(t_0-\epsilon(t_0),t_0+\epsilon(t_0))$, $U(t):=(u_{ij}(t))_{i,j=1}^n$ is a unitary matrix.
			\item[(iii)] for all  $t\in(t_0-\epsilon(t_0),t_0+\epsilon(t_0))$, $U(t)^*(A+tB)U(t)=\text{diag}(\lambda_1(t),\dots,\lambda_n(t)).$ 
		\end{itemize}
	\end{thm}
	The above theorem has been studied extensively. We refer \cite{Ka13}, \cite{ReBe69} and \cite{Ba85} for more on related results on this topic. 
	
	\subsection{Some differentiability criterion:}
	We state and prove the following lemma which will be used later. For any $x\in\mathbb R,$ denote $[x]$ to be the greatest integer such that $[x]\leq x<[x]+1.$ In the following lemma we say a function is zero time differentiable if and only if it is not differentiable at all.
	\begin{lma}\label{difflm}
		Let $p,q$ be two positive real numbers and let $f(t)=(1+|t|^q)^p,\, t\in\R$. Then, at $0$, $f$ is
		infinitely many times differentiable if $q\in2\mathbb{N}$,
			$q -1$  times differentiable if  $q\in2\mathbb{N}-1$, and $[\,q\,]$ times differentiable if  $q\in(0,\infty)\setminus\mathbb{N}$. Moreover, for any two distinct real numbers $q_1,q_2\in(0,1]$, $(1+|t|^{q_1})^p-a|t|^{q_2}$ is not differentiable at $0$, where $a$ is any real constant.
	\end{lma}
	\begin{proof}
		The proof of the first part may be safely left to the reader. For the second part, the binomial series expansion of $(1+|t|^{q_1})^p$ in a small neighbourhood $\Omega$ of $0$ is given by
		\begin{align}\label{l0}
			(1+|t|^{q_1})^p=\sum_{k=0}^{\infty} \left(\begin{matrix}
				p\\
				k
			\end{matrix}\right)|t|^{kq_1},\quad t\in \Omega, \text{ where }\left(\begin{matrix}
				p\\
				k
			\end{matrix}\right)=\frac{p(p-1)\cdots(p-k+1)}{k!}.
		\end{align}
		Let $\bm{m}$ be the least positive integer such that $\bm mq_1> 1$. Since the series in \eqref{l0} is absolutely convergent in $\Omega$, we may write
		\begin{align*}
			(1+|t|^{q_1})^p=\sum_{k=0}^{\bm{m}-1} \left(\begin{matrix}
				p\\
				k
			\end{matrix}\right)|t|^{kq_1}+\mathcal{O}(|t|^{\bm mq_1}),\quad t\in \Omega.
		\end{align*}
		Therefore,
		\begin{align}\label{l1}
			(1+|t|^{q_1})^p-a|t|^{q_2}=\left\{\sum_{k=0}^{\bm{m}-1} \left(\begin{matrix}
				p\\
				k
			\end{matrix}\right)|t|^{kq_1}-a|t|^{q_2}\right\}+\mathcal{O}(|t|^{\bm mq_1}), \quad t\in \Omega.
		\end{align}
		Note that $\left\{\sum\limits_{k=0}^{\bm{m}-1} \left(\begin{matrix}
			p\\
			k
		\end{matrix}\right)|t|^{kq_1}-a|t|^{q_2}\right\}$ is not differentiable at $0$ but $\mathcal{O}(|t|^{\bm mq_1})$ is differentiable at $0$. Hence from \eqref{l1}, we conclude that $(1+|t|^{q_1})^p-a|t|^{q_2}$ is not differentiable. This completes the proof of the lemma.
	\end{proof}

	\section{Isometric Embeddability $\ell_q^m(\mathbb K)\to\ell_p^n(\mathbb K)$ as commutative quasi-Banach space}\label{comiso}
	In this section we study when $\ell_q^m(\mathbb K)$ embeds isometrically into $\ell_p^n(\mathbb K)$ as quasi-Banach spaces.
	\begin{proof}[ Proof of Theorem \ref{main1}]\label{prof thm 1.2}
		On the contrary suppose there exists an isometric embedding of $\ell_q^m(\mathbb{K})$ into $\ell_p^n(\mathbb{K})$. Note that as $\ell_q^2(\mathbb{K})$ is isometrically embedded into $\ell_q^m(\mathbb{K})$, $\ell_q^2(\mathbb{K})$ is also isometrically embedded into $\ell_p^n(\mathbb{K})$. Let $T:\ell_q^2(\mathbb{K})\to\ell_p^n(\mathbb{K})$ be the linear map that embeds $\ell_q^2(\mathbb{K})$ into $\ell_p^n(\mathbb{K})$ isometrically. Let $T(1,0)=\bm{a}:=(a_1,a_2,\ldots,a_n)$ and $T(0,1)=\bm{b}:=(b_1,b_2,\dots,b_n)$. Since $T$ is an isometry, we have
		\begin{align}
			\label{f1}&\|\bm{a}\|_p=\|\bm{b}\|_p=1, \\
			\label{f2}&(1+|t|^q)^{p/q}=\sum_{k=1}^{n}|a_k+tb_k|^p,~~ \forall~t\in\R, \text{ and }\\
			\label{f3}&(1+|t|^q)^{p/q}= \sum_{k=1}^{n}|b_k+ta_k|^p,  ~~ \forall~t\in\R.
		\end{align}
		Let $\bm{m}$ be the least positive integer such that ${\bm m}q> 1$. If all $a_k, 1\leq k\leq n$ are non-zero, then the right hand side of \eqref{f2} is real-analytic in a small neighborhood $\Omega\subset \R$ of $0$, so the left hand side of \eqref{f2} is also real-analytic in $\Omega$. Therefore by Lemma \ref{difflm}, from \eqref{f2} we conclude $q\in 2\mathbb{N}$. So when $q\in (0,\infty)\setminus 2\mathbb{N}$, then this situation never exists. Hence assume $q\in 2 \mathbb{N}$ and let $\mathbb{K}=\R$. Then from \eqref{f1} and \eqref{f2}, we have
		\begin{align*}
			1+\frac{p}{q}\,|t|^{q}+\mathcal{O}(|t|^{2q})=1+t\ddt\Bigg|_{t=0}\|\bm{a}+t\bm{b}\|_p^p+\frac{t^2}{2!} \at{\frac{d^2}{dt^2}}{t=0}\norm{\bm{a}+t\bm{b}}_p^p+\mathcal{O}(|t|^3).
		\end{align*}
	Consequently,
		\begin{align}\label{f4}
			 \frac{p}{q}\,|t|^{q-1}+\mathcal{O}(|t|^{2q-1})=\ddt\Big|_{t=0}\|\bm{a}+t\bm{b}\|_p^p+\frac{t}{2!} \at{\frac{d^2}{dt^2}}{t=0}\norm{\bm{a}+t\bm{b}}_p^p+\mathcal{O}(|t|^2),\quad t\in\Omega.
		\end{align}
		Since $1=\|\bm{a}\|_p^p\leq \|\bm{a}+t\bm{b}\|_p^p$, the function $t\mapsto\|\bm{a}+t\bm{b}\|_p^p$ attains its minimum value at $t=0$. Hence, we must have $\ddt\Big|_{t=0}\|\bm{a}+t\bm{b}\|_p^p=0$. Therefore, from \eqref{f4}, we have 
		\begin{align}\label{f5}
			\frac{p}{q}\,|t|^{q-2}+\mathcal{O}(|t|^{2q-2})={\frac{1}{2!}} \at{\frac{d^2}{dt^2}}{t=0}\norm{\bm{a}+t\bm{b}}_p^p+\mathcal{O}(|t|),\quad t\in\Omega.
		\end{align}
		Taking limit $t\to 0$ on both side of \eqref{f5} we have
		\begin{align}\label{afact}
			\frac{p}{q}\,\lim_{t\to 0}|t|^{q-2}=\frac{1}{2!} \at{\frac{d^2}{dt^2}}{t=0}\norm{\bm{a}+t\bm{b}}_p^p=\frac{p(p-1)}{2!}\sum_{k=1}^{n}|a_k|^{p-2}|b_k|^2<0,
		\end{align}
		which is not possible as $\frac{p}{q}>0$. By similar argument, from $\eqref{f3},$ one can conclude that not all $b_k,  1\leq k\leq n$ are zero. Therefore for the scalar field $\mathbb{K}=\R$ and $q\in(0,\infty)$ or $\mathbb{K}=\mathbb{C}$ and $q\in (0,\infty)\setminus 2\mathbb{N}$, without loss of generality, we may assume that $a_1=\cdots=a_l=0$ for some $l\in\{1,2,\ldots, n-1\}$ and $a_k\neq 0$ for $l+1\leq k\leq n.$ Now if $b_k=0$ for all $k\in\{1,2,\ldots,l\}$, then we can replace $\bm{a}$ and $\bm{b}$ by $(a_{l+1},a_{l+2},\ldots, a_{n})$ and $(b_{l+1},b_{l+2},\ldots,b_{n})$ respectively in \eqref{f1},\eqref{f2} and \eqref{f3}, and by the same analysis as above we can ensure that this case never exists. Note that if $b_k=0$ for $l+1\leq k\leq n$, then for each real number $t$, from \eqref{f2} we have $(1+|t|^q)^{p/q}=1+|t|^p$, which is obviously not true. Therefore, for some $k$ with $1\leq k\leq l,$ $b_k$ is non-zero and for some $s$ with $l+1\leq s\leq n,$ $b_s$ is also non-zero.  Thus, for the scalar field $\mathbb{K}=\R$ and $q\in(0,\infty)$ or $\mathbb{K}=\mathbb{C}$ and $q\in (0,\infty)\setminus 2\mathbb{N}$, from \eqref{f2} we get
		\begin{align}\label{f6}
			(1+|t|^q)^{p/q}-\left(\sum_{k=1}^{l}|b_k|^p\right)|t|^p= \sum_{k=l+1}^{n}|a_k+tb_k|^p,\quad t\in\R.
		\end{align}
		Note that the right hand side of \eqref{f6} is real-analytic in a small neighborhood of $0$, but by Lemma~\ref{difflm}, the left hand side is not differentiable at $0$, which leads to a contradiction. This completes the proof of the theorem.
	\end{proof}
\begin{rmrk}\label{ramk3.1}
		If we consider $\bm{a},\bm{b}\in \ell_p^n(\mathbb{C})$ with $a_k\neq0, 1\leq k \leq n$, such that \eqref{f1} holds. Then it may happen that $\at{\frac{d^2}{dt^2}}{t=0}\norm{\bm{a}+t\bm{b}}_p^p=0$. Indeed, let $\bm a=n^{-\frac{1}{p}}(1,1,\ldots,1)$, $\bm b=\frac{n^{-\frac{1}{p}}}{\sqrt{2-p}}\big(1+i\sqrt{1-p}, \ldots, 1+i\sqrt{1-p}\big)$. Then note that $\|\bm a\|_p=\|\bm b\|_p=1$, and 
		\[\at{\frac{d^2}{dt^2}}{t=0}\norm{\bm{a}+t\bm{b}}_p^p=\at{\frac{d^2}{dt^2}}{t=0}\left(t^2+\frac{2}{\sqrt{2-p}}t+1\right)^{\frac{p}{2}}=0.\]
		So in this case from \eqref{f5} we can't conclude anything about isometric embeddability of $\ell_q^m(\mathbb{C})$ into $\ell_p^n(\mathbb{C})$ as the above proof breaks down at \ref{afact}.
	\end{rmrk}
	\section{Isometric Embeddability $S_q^m\to S_p^n$ as non-commutative quasi-Banach spaces}\label{ncomiso}
	In this section we study when $S_q^m$ embeds isometrically into $S_p^n$ as quasi-Banach spaces. We also give applications of our result to non-commutative $L_p$-spaces. We begin recalling the concepts of operator derivatives and multiple operator integrals which were used in \cite{JFA22}.
	\subsection{ Operator derivatives in $\ell_2^m$ in terms of multiple operator integrals.}
\begin{dfn}
		Let $M_0,M_1,\ldots,M_n$ be self-adjoint operators on $\ell_2^m$, let $\mathfrak{e}^{(j)}=\{\mathfrak{e}_i^{(j)}\}_{i=1}^{m}$ be an orthonormal basis of eigenvectors of $M_j$, and let $\bm{\lambda}^{(j)}=(\lambda_i^{(j)})_{i=1}^{m}$ be the corresponding $m$-tuple of eigenvalues for $j=0,1,\ldots,n$. Let $\Phi:\R^{n+1}\to\mathbb{C}$. We define $$T^{M_0,M_1,\ldots,M_n}_\Phi: \underbrace{\blm\times\blm\times\cdots\times\blm}_{\text{$n$-times}} \longrightarrow \blm$$ by 
	\begin{align}\label{01}
		T^{M_0,M_1,\ldots,M_n}_\Phi(N_1,N_2,\ldots,N_n)=\sum_{i_0,i_1,\ldots,i_n=1}^{m}\Phi(\lambda_{i_0}^{(0)},\lambda_{i_1}^{(1)},\ldots,\lambda_{i_n}^{(n)})P_{\mathfrak{e}_{i_0}^{(0)}}N_1P_{\mathfrak{e}_{i_1}^{(1)}}N_2\cdots N_nP_{\mathfrak{e}_{i_n}^{(n)}},
	\end{align}
	for any $n$-tuple $(N_1,N_2,\ldots,N_n)$ in $\blm$, where $P_{\mathfrak{e}_{i_k}^{(k)}}$ is the orthogonal projection of $\ell_2^m$ onto the subspace spanned by the vectors $\{\mathfrak{e}_{i_k}^{(k)}\}$. The operator $T^{M_0,M_1,\ldots,M_n}_\Phi$ is a discrete version of a multiple operator integral and the function $\Phi$ is called the symbol of the operator $T^{M_0,M_1,\ldots,M_n}_\Phi$.
	
\end{dfn}
	For additional information on discrete multiple operator integral and related topics, we refer \cite[Chapter 4]{SkTo19}.
	
	Throughout this paper we denote by $C^k(\R)$ to be the space of all $k$-times differentiable functions which also has continuous $k$-th derivative, where $k\in\mathbb N\cup\{0\}$. Recall that the divided difference of order $r$ is an operation on a function $f\in C^k(\R)$ defined recursively as follows:
	\begin{align*}
		f^{[0]}(\lambda)&:=f(\lambda),\\
		f^{[r]}(\lambda_0,\lambda_1,\ldots,\lambda_r)&:=\begin{cases*}
			\frac{f^{[r-1]}(\lambda_0,\lambda_1,\ldots,\lambda_{r-2},\lambda_r)-f^{[r-1]}(\lambda_0,\lambda_1,\ldots,\lambda_{r-2},\lambda_{r-1})}{\lambda_r-\lambda_{r-1}} \quad \text{if}\quad \lambda_r\neq\lambda_{r-1},\\
			\at{\frac{\partial}{\partial \lambda}}{\lambda=\lambda_{r-1}}f^{[r-1]}(\lambda_0,\lambda_1,\ldots,\lambda_{r-2},\lambda)\quad \text{if}\quad \lambda_r=\lambda_{r-1}.
		\end{cases*}
	\end{align*}
	We recall the formulae for operator derivatives in terms of multiple operator integrals.
	\begin{thm}\cite[Theorem 5.3.2]{SkTo19}\label{diffor}
		Let $M,N\in M_n^{sa}$ and $f\in C^r(\R), r\in\Nat$. Then the function $\mathbb{R} \ni t\mapsto f(M+tN)$ is $r$-times differentiable in the operator norm and the $r$-th order derivative is given by the formula
		\begin{align}
			\label{eq2}&\frac{1}{r!}\at{\frac{d^r}{dt^r}}{t=0}f\big(M+tN \big)=T^{M,M,\ldots,M}_{f^{[r]}}\underbrace{(N,N,\ldots,N)}_{r\text{-times}}
		\end{align}
		and hence 
		\begin{align}
			\label{eq3}&\frac{1}{r!}\at{\frac{d^r}{dt^r}}{t=0}\Tr\big(f\big(M+tN \big)\big)=\Tr\bigg(T^{M,M,\ldots,M}_{f^{[r]}}\underbrace{(N,N,\ldots,N)}_{r\text{-times}}\bigg).
		\end{align}
	\end{thm}
	
	By applying Theorem \eqref {diffor}, one can get the following first-order derivative formula; however, we refer to this formula here because this formula can also be proved by a simple application of the Kato-Rellich theorem.
	\begin{ppsn}
		Let $A,B\in M_n^{sa}$ with $A$ is invertible and $0<p\leq 1$, then 
		\begin{align*}
			\at{\frac{d}{dt}}{t=0}\norm{A+tB}_p^p=p\Tr\big(B|A|^{p-1}\sgn(A)\big).
		\end{align*}
	\end{ppsn}
	\begin{proof}
		The result follows from the same lines of proof of \cite[Proposition 2.6]{JFA22}.
	\end{proof}
	
	\subsection{Reduction to self-adjoint case:} We now state the following lemma, which reduces our problem to the self-adjoint case. The following result is a straightforward generalization of \cite[Lemma 2.7]{JFA22}. Hence we omit the proof of the following lemma. 
	\begin{lma}\label{reduct}
		Let $0< q \leq \infty$, $0<p<1$. Let $T:\ell_q^2\to S_p(\ell_2^n)$ with $T(e_1)=A,$ $T(e_2)=B$ be an isometry. Then \[T_{new}:\ell_q^2\to S_p(\ell_2^n\oplus\ell_2^n)\text{ defined by } T_{new}(z,w)\colon=zA_{new}+wB_{new}\] is again an isometry, where \[A_{new}(\zeta_1\oplus\zeta_2)\colon=2^{-\frac{1}{p}}\Big(A\zeta_2\oplus A^*\zeta_1\Big)~ \text{and}~ B_{new}(\zeta_1\oplus\zeta_2)\colon=2^{-\frac{1}{p}}\Big(B\zeta_2\oplus B^*\zeta_1\Big).\]
	\end{lma}
	\begin{dfn}\label{hp}
		Let $0< q \leq \infty$, $0<p\leq \infty$, and $p\neq q$. Let $A,B\in \mathcal{B}(\ell_2^n).$ We say $(A,B)$ has $(\mathbf{I}_{q,p})$ if
		\begin{itemize}
			\item[$\bullet$] with respect to the standard basis, $A$ is a diagonal matrix with all entries real and $B$ is self-adjoint,
			\item[$\bullet$] $\norm{A}_p=\norm{B}_p=1$,
			\item[$\bullet$] $\norm{A+tB}_p=\|(1,t)\|_q$, for all $t\in\R$.
		\end{itemize}
	\end{dfn}
	
	We prove the following lemma.
	\begin{lma}\label{l}
		Let $0< q \leq \infty$, $0<p<1$, $p\neq q$. Assume that $A, B\in L_p(\mathcal M)$ are self-adjoint elements with $\|A\|_p=\|B\|_p=1$ and $\norm{A+tB}_p=\|(1,t)\|_q$, for all $t\in\R$. Then we must have $AB\neq 0.$
	\end{lma}
	\begin{proof}
		Suppose $AB=0$. Then $A$ and $B$ have disjoint supports. It follows that $\|A+tB\|_p^p=\|A\|_p^p+|t|^p\|B\|_p^p$. Note that $\|A\|_p=\|B\|_p=1$. Thus $\|A+tB\|_p^p=1+|t|^p.$  But for each $t\in \R$, the identity $(1+|t|^q)^{1/q}=(1+|t|^p)^{1/p}$ is not true, which leads to a contradiction.
	\end{proof} 
	\par Due to the above, we have the following necessary condition for existence of an isometric embedding of $S_q^m$ into $S_p^n$ for $0< q \leq \infty$, $0<p<1$, $p\neq q$.
	\begin{crlre}\label{neq0}
		Let $0< q \leq \infty$, $0<p<1$, $p\neq q$ and $(A,B)$ has $(\mathbf{I}_{q,p})$. Then we must have that $AB\neq 0$. 
	\end{crlre}
	The following theorem is one of the main result of this section. One of the main difficulties we face  to prove the non-existence of isometric embedding is the nature of differentiability of the map $t\mapsto \|A+tB\|_p^p$ at $t=0,1$, which is keenly related to the multiplicity of zeros of the real-analytic eigenvalues of $A+tB$ at $t=0$ and $t=1$, where $(A,B)$ has $(\mathbf{I}_{q,p})$. We closely examine the multiplicity of zeros and compute the derivatives using the Kato-Rellich theorem and multiple operator integrals. This helps us to analyze the power-series of various analytic functions and compare them carefully. We need to use \textsf{repeatedly} the Kato-Rellich theorem and compare the coefficients arising in the series expansion. This also leads to various intricate computations which are combinatorial in nature.
	
	\begin{thm}\label{th1}
		Let $0<p\neq q<1$ and $2\leq m\leq n<\infty$. Then there is no isometric embedding of $S_q^m$ into $S_p^n$. 
	\end{thm}
	\begin{proof}
		On the contrary, suppose that there exists an isometric embedding $S_q^m$ into $S_p^n$. Then by Lemma \ref{lem2.2}, we have $p<q$. Therefore, by Lemma \ref{reduct}, we may assume that there exists a diagonal self-adjoint matrix $A$ and a self-adjoint matrix $B$ in $S_p^n$ such that $(A, B)$ has $(\mathbf{I}_{q,p})$. In particular, for any real number $t$, we have the identity
		\begin{align}\label{a1}
			\left(1+|t|^q\right)^{\frac{p}{q}}=\|A+tB\|_p^p.
		\end{align}
		In a small neighborhood $\Omega$ of zero, by Kato-Rellich Theorem \ref{kato}, there exist $n$ real-analytic functions $\lambda_1(\cdot), \lambda_2(\cdot),\ldots, \lambda_n(\cdot)$ on $\Omega$ such that they are complete set of eigenvalues of $A+tB$ for $t\in\Omega$. So for $t\in\Omega $, from \eqref{a1} we have 
		\begin{align}\label{a2}
			\left(1+|t|^q\right)^{\frac{p}{q}}=\sum_{k=1}^{n}|\lambda_k(t)|^p.
		\end{align}
		If  $\lambda_k(0)\neq 0$ for  all $1\leq k\leq n$, that is if $A$ is invertible, then the right hand side of $\eqref{a2}$ is a real-analytic function in $\Omega$ but the left hand side of \eqref{a2} is not analytic in $\Omega$. So $A$ must be a singular matrix. Again from \eqref{a2}, it is clear that not all $\lambda_k(0)$ are zero. So there is a nonempty proper subset $O\subset \{1,2,\ldots,n\}$ such that $\lambda_k(0)=0$ for $k\in O$, and $\lambda_k(0)\neq 0$ for $k\in \{1,2,\ldots,n\}\setminus O$. Since the quasi-norm $\|\cdot\|_p$ is unitary invariant, we may assume that there exists an $l\in\{1,2,\ldots,n-1\}$ such that $\lambda_k(0)=0$ for $1\leq k\leq l$, and $\lambda_k(0)\neq 0$ for $l+1\leq k\leq n$. Note that if some $\lambda_i(t)$ are identically zero, we can neglect those as they have no contribution in $\|A+tB\|_p$. Therefore, we may also assume that all $\lambda_k(t)$ are non-zero analytic functions in $\Omega$. Let $m_k$ be the multiplicity of zero at $0$ of the function $\lambda_k(t)$ for $1\leq k\leq l$. Then
		\begin{align}\label{a3}
			\lambda_k(t)=t^{m_k}\,\mu_k(t),~~~ \mu_k(0)\neq 0,~~ 1\leq k\leq l,
		\end{align}
		where $\mu_k(\cdot)$ are analytic in $\Omega$. For further progression, let us define
		\begin{align}\label{a4}
			\Psi(t)=\sum_{k=l+1}^{n}|\lambda_k(t)|^p-1, \quad t\in\Omega.
		\end{align} 
		It is clear that $\Psi(\cdot)$ is analytic in $\Omega$ and $\Psi(0)=0$. Then $\Psi(t)=t\Psi_{1}(t)$, where $\Psi_{1}(\cdot)$ is a analytic function in $\Omega$. Applying \eqref{a3} and \eqref{a4} in \eqref{a2}, we have
		\begin{align}\label{a5}
			\left(1+|t|^q\right)^{\frac{p}{q}}=\sum_{k=1}^{l}|t|^{m_k p}|\,\mu_k(t)|^p+t\Psi_{1}(t)+1, \quad t\in\Omega.
		\end{align}
		Now we will reach our goal through case-by-case analysis.\\[5pt]
		
		{\bf Case I:} Let $m_kp>1$ for $1\leq k\leq l$. Then the right hand side of \eqref{a5} is differentiable at $0$ but the left hand side of \eqref{a5} is not differentiable at $0$. So this case is not possible.\\[5pt]
		
		{\bf Case II:} Suppose at least one of $m_kp=1$ for $1\leq k\leq l$. Without loss of generality, we assume that $m_1p=1.$ Note that as $\mu_k(0)\neq 0$,  $|\mu_k(t)|^p$ are analytic in $\Omega$ for $1\leq k\leq l$. So we can express $|\mu_k(t)|^p$ as
		\begin{align}\label{a6}
			|\mu_k(t)|^p=|\mu_k(0)|^p+ t\,\xi_k(t), ~~1\leq k\leq l, \quad t\in\Omega,
		\end{align}
		where $\xi_k(\cdot)$ are some analytic function in the small neighborhood $\Omega$ of $0$. Let $\bm{m}$ be the least positive integer such that $\bm mq> 1$.
		
		{\bf Sub-case I:}
		 Suppose $l\geq 3$, $m_1p=m_2p=\cdots=m_jp=1$, $m_{j+1}p, \ldots, m_ip<1$, and $m_{i+1}p,\ldots,m_lp>1$ for some $j\in\{1,2,\ldots, l-2\}$, and $i\in\{j+1\ldots, l-1\}$. Combining \eqref{a5} and \eqref{a6}, we have 
		\begin{align}\label{a7}
				\nonumber&1+\alpha_1\,|t|^{q}+\alpha_2\,|t|^{2q}+\cdots+\alpha_{\bm{m}-1}\,|t|^{(\bm{m}-1)q}+\mathcal{O}(|t|^{\bm mq})
				=1+\sum_{k=1}^{j}|t|\,|\mu_k(0)|^p\\
				&+\sum_{k=j+1}^{i} |t|^{m_kp}\,|\mu_k(0)|^p
				+\sum_{k=1}^{i} |t|^{m_kp}t\xi_k(t)+\sum_{k=i+1}^{l} |t|^{m_kp}\,|\mu_k(t)|^p+t\Psi_{1}(t), \quad t\in\Omega,
		\end{align}
		where $\alpha_k$ are the coefficient of $|t|^{kq}$ in the binomial expansion of $(1+|t|^q)^{p/q}$ for $|t|<1$. That is $\alpha_k=\left(\begin{matrix}
			p/q\\
			k
		\end{matrix}\right)=\frac{p/q(p/q-1)\cdots(p/q-k+1)}{k!}$, $k\in \mathbb{N}$.\\
		From \eqref{a7} we have
		\begin{align}\label{a8}
				\nonumber&\Big\{\alpha_1\,|t|^{q}+\alpha_2\,|t|^{2q}+\cdots+\alpha_{\bm{m}-1}\,|t|^{(\bm{m}-1)q}\Big\}-\Big\{\sum_{k=1}^{j}|t|\,|\mu_k(0)|^p + \sum_{k=j+1}^{i} |t|^{m_kp}\,|\mu_k(0)|^p           \Big\}\\
				=& \sum_{k=1}^{i} |t|^{m_kp}t\xi_k(t) + \sum_{k=i+1}^{l} |t|^{m_kp}\,|\mu_k(t)|^p+\mathcal{O}(|t|^{\bm mq})+t\Psi_{1}(t), \quad t\in\Omega.
		\end{align}
		In this case, the right hand side of the above equation \eqref{a8} is differentiable at $0$. Since $p<q$, so $\alpha_{2}<0$. Therefore by the similar lines of argument as given in the proof of Lemma \ref{difflm}, the left hand side of \eqref{a8} is not differentiable at $0$. Hence this case does not exist.\\
		
		{\bf Sub-case II:} Suppose $l\geq 2$, $m_1p=\cdots=m_jp=1$, and $ m_{j+1}p,\ldots,m_lp<1$ for some $j\in\{1,\ldots,l-1\}$. Then for $t\in \Omega$, from \eqref{a7},  we get
			\begin{align}\label{a9}
				\nonumber&\alpha_1\,|t|^{q}+\cdots+\alpha_{\bm{m}-1}\,|t|^{(\bm{m}-1)q}
				-\sum_{k=1}^{j}|t|\,|\mu_k(0)|^p-\sum_{k=j+1}^{l} |t|^{m_kp}\,|\mu_k(0)|^p\\
				=&\sum_{k=1}^{l} |t|^{m_kp}t\xi_k(t)+\mathcal{O}(|t|^{\bm mq})+t\Psi_{1}(t).
		\end{align}
		It is clear that the left side of the equation \eqref{a9} is not differentiable at $0$, but the other side is differentiable at $0$. So the equation \eqref{a9} never exists in $\Omega$.
		
		{\bf Sub-case III:} Suppose $l\geq 2$, $m_1p=\cdots=m_jp=1$, $ m_{j+1}p,\ldots,m_lp>1$ for some $j\in\{1,\ldots,l-1\}$. Again from \eqref{a7} we have 
			\begin{align}\label{a10}
				\nonumber&\alpha_1\,|t|^{q}+\cdots+\alpha_{\bm{m}-1}\,|t|^{(\bm{m}-1)q}-\sum_{k=1}^{j}|t|\,|\mu_k(0)|^p\\
				=& \sum_{k=1}^{j}|t|\, t\, \xi_k(t)+\sum_{k=j+1}^{l} |t|^{m_kp}\,|\mu_k(t)|^p+t\Psi_{1}(t)+\mathcal{O}(|t|^{\bm mq}), ~~ t\in\Omega.
		\end{align} By a similar argument as given for the nonexistence of equation \eqref{a9}, we conclude the nonexistence of equation \eqref{a10}.\\
		
		{\bf Sub-case IV:} Let $m_1p=\cdots=m_lp=1$. By similar argument as given in the above sub-case, we conclude that this situation never exists.
		
		\vspace{.1in}
		{\bf Case III:} Let $m_kp<1$ for $1\leq k\leq l$. From \eqref{a7} we have
		\begin{align}\label{a11}
			\alpha_1\,|t|^{q}+\cdots+\alpha_{\bm{m}-1}\,|t|^{(\bm{m}-1)q}
			-\sum_{k=1}^{l} |t|^{m_kp}\,|\mu_k(0)|^p=\sum_{k=1}^{l} |t|^{m_kp}t\xi_k(t)+t\Psi_{1}(t)+\mathcal{O}(|t|^{\bm mq}), ~ t\in\Omega.
		\end{align}
		Note that the right hand side of \eqref{a11} is differentiable at $0$. Therefore the existence of the equation \eqref{a11} implies that the left hand side of \eqref{a11} is also differentiable at $0$. \\
		
		{\bf Sub case I:}
		Let $2q\leq 1$. The necessary condition for differentiability of the left hand side of \eqref{a11} is that there exists a nonempty proper subset $I\subset\{1,2,\ldots,l\}$ such that 
			\begin{enumerate}[(A)]
				\item\label{id} $m_kp=2q,~~ k\in I,$ and $\alpha_2=\sum\limits_{k\in I}|\mu_k(0)|^p.$
			\end{enumerate}
			By Lemma \ref{lem2.2}, we obtain $\frac{p}{q}<1$. Hence $\alpha_2=\frac{1}{2!}\frac{p}{q}\left(\frac{p}{q}-1\right)<0$, but $\sum\limits_{k\in I}|\mu_k(0)|^p>0$. So the above identity \eqref{id} is not true. Hence the equation \eqref{a11} never exists in the small neighborhood $\Omega$, which is a contradiction.\\
		
		{\bf Sub case II:}
		Let $2q> 1$. In this case, the necessary condition for differentiability of the left hand side of \eqref{a11} is that
		\begin{enumerate}[(a)]
			\item\label{idd} $m_kp=q,~~ k\in \{1,2,\ldots,l\}, \text{ and }  \alpha_1=\sum\limits_{k=1}^{l}|\mu_k(0)|^p$.
		\end{enumerate}
		Suppose the above necessary condition \eqref{idd} hold, then from \eqref{a11} we have,
		\begin{align}\label{a12}
			\alpha_{2}\,|t|^{2q-1}\sgn(t)+\alpha_{3}\,|t|^{3q-1}\sgn(t)+\mathcal{O}(|t|^{4q-1})=\sum_{k=1}^{l} |t|^{q}\xi_k(t)+\Psi_{1}(t), \quad t\in\Omega\setminus\{0\}.
		\end{align}
		Taking limit $t\to 0$ on both side of \eqref{a12}, we have $\Psi_{1}(0)=0$ which implies $\Psi_{1}(t)=t\Psi_{2}(t)$ for some analytic function $\Psi_{2}(t)$ in $\Omega$ of $0$. Let $\xi_k(t)=\xi_k(0)+t\,\eta_k(t), 1\leq k\leq l,$ and $\eta_k$ are some analytic functions in $\Omega$. Then from \eqref{a12} we have
		\begin{align}\label{a13}
			\alpha_{2}\,|t|^{2q-2}+\alpha_{3}\,|t|^{3q-2}-|t|^{q-1}\sgn(t)\sum_{k=1}^{l} \xi_k(0)=\sum_{k=1}^{l}|t|^q \eta_k(t)+\Psi_{2}(t)+\mathcal{O}(|t|^{4q-2}), ~ t\in\Omega_2\setminus\{0\}.
		\end{align}
		Observe that the right hand side of \eqref{a13} is bounded around $0$ but the left hand side is unbounded around $0$, which leads to a contradiction.\\
		
		Hence {\bf Case III} does not exist.\\
		
		{\bf Case IV:} Suppose $l\geq 2$ and there are some $m_kp$ which are strictly less than $1$ and rest are strictly greater than $1$. We may assume  $ m_1p,\ldots,m_ip<1$ and $m_{i+1}p,\ldots,m_lp>1$ for some $i\in\,\{1,2,\ldots,l-1\}$.  From \eqref{a7}, we obtain
			\begin{align}\label{h1}
				\nonumber&\Big\{\alpha_1\,|t|^{q}+\alpha_2\,|t|^{2q}+\cdots+\alpha_{\bm{m}-1}\,|t|^{(\bm{m}-1)q}\Big\}-\sum_{k=1}^{i} |t|^{m_kp}\,|\mu_k(0)|^p\\
				=& \sum_{k=1}^{i} |t|^{m_kp}t\xi_k(t) + \sum_{k=i+1}^{l} |t|^{m_kp}\,|\mu_k(t)|^p+\mathcal{O}(|t|^{\bm mq})+t\Psi_{1}(t), \quad t\in\Omega.
			\end{align}
			Then we again consider the following two sub-cases:\\
			
			{\bf Sub case I:} Let $2q\leq 1$. By the similar lines of argument as given in {\bf Sub case I} of {\bf Case III} of this proof, we conclude that this situation never exists.\\
			
			{\bf Sub Case II: } Let $2q> 1$. In this case, the necessary condition for differentiability of the left hand side of \eqref{a11} is that
			\begin{enumerate}[(1)]
				\item\label{s} $m_kp=q,~~ k\in \{1,2,\ldots,i\}, \text{ and }  \alpha_1=\sum\limits_{k=1}^{i}|\mu_k(0)|^p$.
			\end{enumerate}
			Suppose the above necessary condition \eqref{s} hold, then from \eqref{h1} we have,
			\begin{align}\label{h2}
				\nonumber&\alpha_{2}\,|t|^{2q-1}\sgn(t)+\alpha_{3}\,|t|^{3q-1}\sgn(t)+\mathcal{O}(|t|^{4q-1})\\
				=& \sum_{k=1}^{i} |t|^{q}\xi_k(t) + \sum_{k=i+1}^{l} |t|^{m_kp-1}\,\sgn(t)|\mu_k(t)|^p+\Psi_{1}(t), \quad t\in\Omega.
			\end{align}
			Since $2q>1$ and $m_kp>1$ for all $i+1\leq k\leq l$, then taking limit $t\to 0$ on both sides of \eqref{h2}, we have $\Psi_1(0)=0$. Let $\Psi_{1}(t)=t\Psi_{2}(t)$. Let $\xi_k(t)=\xi_k(0)+t\,\eta_k(t), 1\leq k\leq l,$ and $\eta_k$ are some analytic functions in $\Omega$. Then from \eqref{h2}, we get
			\begin{align}\label{h3}
				\nonumber&\alpha_{2}\,|t|^{2q-2}+\alpha_{3}\,|t|^{3q-2}-\sgn{(t)}\,|t|^{q-1}\,\sum_{k=1}^{i} \xi_k(0)-\sum_{k=i+1}^{l} |t|^{m_kp-2}\,|\mu_k(0)|^p\\
				=& \sum_{k=1}^{i} |t|^{q}\eta_k(t) + \sum_{k=i+1}^{l} |t|^{m_kp-2}\,t\,\xi_k(t)+\Psi_{2}(t) +\mathcal{O}(|t|^{4q-2}), \quad t\in\Omega\setminus\{0\}.
			\end{align}
			Note that as $\alpha_{2}<0$, the left hand side of \eqref{h3} is unbounded near $0$ but the right hand side of \eqref{h3} is bounded near $0$ and tends to $\Psi_2(0)$ as $t\to 0$, which implies that this case never exists.

		Therefore from all the above cases we conclude that there does not exist any isometric embedding of $S_q^m$ into $S_p^n$ for $0<q\neq p<1$. This completes the proof of the theorem.
	\end{proof}

	\begin{thm}\label{th2}
		Let $(q,p)\in\left(1,2\right)\times\left(0,1\right) $  and $2\leq m\leq n<\infty$. Then there does not exist any isometric embedding of $S_q^m$ into $S_p^n$ with $p\neq q$.
	\end{thm}
	\begin{proof}
		Suppose that there is an isometric embedding of $S_q^m$ into $S_p^n$.  Then there exist $A,B\in S_p^n$ such that $(A,B)$ has $(\mathbf{I}_{q,p})$. Here we use the same notations used in the proof of Theorem~\ref{th1}. So the existence of isometric embedding implies that the identity \eqref{a2} is true in a small neighbourhood $\Omega$ of $0$. Note that if $A$ is invertible then the right hand side of \eqref{a2} is real-analytic but the left hand side of \eqref{a2} is not real-analytic in $\Omega$. So the matrix $A$ is not invertible. Therefore by the similar argument as given in the proof of Theorem~\eqref{th1}, we can deduce \eqref{a5} from the identity \eqref{a2}. In other words, we have the following.
		\begin{align}\label{b1}
			\left(1+|t|^q\right)^{\frac{p}{q}}=\sum_{k=1}^{l}|t|^{m_k p}|\,\mu_k(t)|^p+t\Psi_{1}(t)+1,\quad t\in\Omega.
		\end{align}
		Note that the left hand side of \eqref{b1} is one time differentiable at $0$, and one time differentiability of \eqref{b1} implies that some of $m_kp$ are in $(1,2)$ and rest belong to $[2,\infty)$. Since $\|\cdot\|_p$ is unitary invariant, we may assume that there exists $i\in\{1,2,\ldots,l-1\}$ such that 
		\begin{align*}
			1<m_kp<2 \text{ for } 1\leq k\leq i, \text{ and } 2\leq m_kp<\infty \text{ for } i+1\leq k\leq l.
		\end{align*} 
		A simple calculation shows that $\Psi_{1}(0)=0$. Let $|\mu_k(t)|^p=|\mu_k(0)|^p+ t\xi_k(t)$, $1\leq k\leq l$, and $\Psi_{1}(t)=t\Psi_{2}(t)$, where $\xi_k(\cdot)$ and $\Psi_{2}(\cdot)$ are some analytic function in $\Omega$. Then from \eqref{b1}, we have 
		\begin{align}\label{b2}
			\nonumber&1+\alpha_1\,|t|^{q}+\alpha_2\,|t|^{2q}+\mathcal{O}(|t|^{3q})\\
			\nonumber=&\sum_{k=1}^{i} |t|^{m_kp}\,|\mu_k(0)|^p+\sum_{k=1}^{i} |t|^{m_kp}t\xi_k(t)+\sum_{k=i+1}^{l} |t|^{m_kp}\,|\mu_k(t)|^p+t^2\Psi_{2}(t)+1,\quad t\in\Omega_1,\\
			\nonumber\implies&\alpha_1\,|t|^{q-2}-\sum_{k=1}^{i} |t|^{m_kp-2}\,|\mu_k(0)|^p\\
			=&\sum_{k=1}^{i} |t|^{m_kp-2}t\xi_k(t)+\sum_{k=i+1}^{l} |t|^{m_kp-2}\,|\mu_k(t)|^p+\Psi_{2}(t) -\alpha_2\,|t|^{2q-2}+\mathcal{O}(|t|^{3q-2}),\quad t\in\Omega\setminus\{0\}.
		\end{align}
		Observe that if $m_kp\neq q$ for some $k\in\{1,\ldots,i\}$, then the left hand side of above equation \eqref{b2} becomes unbounded near $0$ but the right hand side of $\eqref{b2}$ is bounded and going to $\Psi_2(0)$ as $t\to \infty$. So $m_kp=q$ for all $k\in\{1,\ldots,i\}$. Again by the same reason we have $\alpha_1=\sum\limits_{k=1}^{i}|\mu_k(0)|^p$. Therefore from \eqref{b2} we have
		
		\begin{align}\label{b3}
			\alpha_2\,|t|^{2q-2}=\sum_{k=1}^{i} |t|^{q-2}t\xi_k(t)+\sum_{k=i+1}^{l} |t|^{m_kp-2}\,|\mu_k(t)|^p+\Psi_{2}(t) +\mathcal{O}(|t|^{3q-2}),\quad t\in\Omega.
		\end{align}
		Now if there is a nonempty set $I\subseteq\{i+1,\ldots,l\}$ such that $m_kp=2$ for $k\in I$, then we can rewrite \eqref{b3} as
		\begin{align*}
			\alpha_2\,|t|^{2q-2}=\sum_{k=1}^{i} |t|^{q-2}t\xi_k(t)+\sum_{k\in \{i+1,\ldots,l\}\setminus I} |t|^{m_kp-2}\,|\mu_k(t)|^p+\tilde{\Psi}_{2}(t) +\mathcal{O}(|t|^{3q-2}),\quad t\in\Omega,
		\end{align*}
		where $\tilde{\Psi}_{2}(t)=\Psi_{2}(t)+\sum_{k\in  I} |\mu_k(t)|^p,\, t\in\Omega$. Clearly $\tilde{\Psi}_2(\cdot)$ is analytic in $\Omega$. Therefore in \eqref{b3}, without loss of generality, we may assume that $m_kp>2$ for all $i+1\leq k\leq l$ as $\Psi_2$ can be replaced by $\tilde{\Psi}_2$ and other term can be modified accordingly.
		Taking limit $t\to 0$ on both sides of \eqref{b3}, we conclude that $\Psi_{2}(0)=0$. Let $\xi_k(t)=\xi_k(0)+t\eta_k(t), 1\leq k\leq l$, and $\Psi_{2}(t)=t\Psi_{3}(t)$, where $\eta_k(\cdot)$, and $\Psi_{3}(\cdot)$ are some analytic functions in $\Omega$. Then from \eqref{b3}, we have 
		\begin{align}\label{b4}
			\nonumber&\alpha_2\,|t|^{2q-2}-\sum_{k=1}^{i} |t|^{q-2}\, t \, \xi_k(0)\\
			=&\sum_{k=i+1}^{l} |t|^{m_kp-2}\,|\mu_k(0)|^p+\sum_{k=1}^{i} |t|^{q-2}\, t^2\, \eta_k(t)+\sum_{k=i+1}^{l} |t|^{m_kp-2}t \,\xi_k(t)+t\Psi_{3}(t) +\mathcal{O}(|t|^{3q-2}),\quad t\in\Omega.
		\end{align}
		\textbf{Case I:} Let $1<q\leq \frac{3}{2}$, then $2q-2\leq 1$. Therefore $|t|^{2q-1}$ is not differentiable at $0$. Now if $m_kp>3$ for all $k\in\{i+1,\ldots,l\}$ then the right hand side of \eqref{b4} is differentiable at $0$ but the left hand side  of \eqref{b4} is not differentiable at $0$. So without loss of any generality, we may assume that there exists $j\in\{i+1,\ldots,l-1\}$ such that $m_kp\leq 3$ for $i+1\leq k\leq j$ and $m_kp>3$ for $j+1\leq k\leq l$. Then from \eqref{b4}, we have
		\begin{align}\label{b5}
			\nonumber&\alpha_2\,|t|^{2q-2}-|t|^{q-2}\, t\sum_{k=1}^{i} \xi_k(0)-\sum_{k=i+1}^{j} |t|^{m_kp-2}\,|\mu_k(0)|^p\\
			=&\sum_{k=j+1}^{l} |t|^{m_kp-2}\,|\mu_k(0)|^p+\sum_{k=1}^{i} |t|^{q-2}\, t^2\, \eta_k(t)+\sum_{k=i+1}^{l} |t|^{m_kp-2}t \,\xi_k(t)+t\Psi_{3}(t) +\mathcal{O}(|t|^{3q-2}),~ t\in\Omega.
		\end{align}
		Note that the right hand side of \eqref{b5} is differentiable at $0,$ which forces us to conclude that there exists a non-empty subset $J\subseteq \{i+1,\ldots,j\}$ such that 
		$m_kp=2q$ for $k\in J$ and $\alpha_{2}=\sum\limits_{k\in J} |\mu_k(0)|^p$. But $\alpha_{2}<0$ and $\sum\limits_{k\in J} |\mu_k(0)|^p>0$, which leads to a contradiction.\\
		
		\textbf{Case II:} Let $\frac{3}{2}<q<2$. Then from \eqref{b4} we have 
		\begin{align}\label{b6}
			\nonumber&\alpha_2\,\frac{|t|^{2q}}{t^3}-\left(\sum_{k=1}^{i} \xi_k(0)\right)|t|^{q-2}\,
			=\,\sum_{k=i+1}^{l} \frac{|t|^{m_kp}}{t^3}\,|\mu_k(0)|^p+\left(\sum_{k=1}^{i} \eta_k(t)\right) \frac{|t|^{q}}{t}\\
			&\hspace{3cm}+\sum_{k=i+1}^{l} |t|^{m_kp-2} \,\xi_k(t)+\Psi_{3}(t) +\mathcal{O}(|t|^{3q-3}), \quad t\in\Omega\setminus\{0\}.
		\end{align}
		Let $A=\sum\limits_{k=1}^{i} \xi_k(0)$, and $\eta_k(t)=\eta_k(0)+t\,\gamma_k(t),\, 1\leq k\leq l$, for some analytic function $\gamma_k(\cdot)$ in $\Omega$. Now we rewrite the above equation \eqref{b6} as
		\begin{align}\label{b7}
			\nonumber&\alpha_2\,\frac{|t|^{2q}}{t^3}-A|t|^{q-2}-\left(\sum_{k=1}^{i} \eta_k(0)\right) \frac{|t|^{q}}{t}\,=\,\sum_{k=i+1}^{l} \frac{|t|^{m_kp}}{t^3}\,|\mu_k(0)|^p + \left(\sum_{k=1}^{i} \gamma_k(t)\right) |t|^{q}\\
			&\hspace{3cm}+\sum_{k=i+1}^{l} |t|^{m_kp-2} \,\xi_k(t)+\Psi_{3}(t) +\mathcal{O}(|t|^{3q-3}), \quad t\in\Omega\setminus\{0\}.
		\end{align}
		If $m_kp\geq 4$ for $i+1\leq k\leq l$, then the right hand side of \eqref{b7} is differentiable but the left hand side of \eqref{b7} is not differentiable at $0$. So we may assume that there exist some $j\in\{i+1,\ldots,l-1\}$ such that $2\leq m_kp<4$ for $i+1\leq k\leq j,$ and $4\leq m_kp<\infty$ for $j+1\leq k\leq l.$ Then from \eqref{b7} we have
		\begin{align}\label{b8}
			\nonumber&\alpha_2\,\frac{|t|^{2q}}{t^3}-A|t|^{q-2}-\left(\sum_{k=1}^{i} \eta_k(0)\right) \frac{|t|^{q}}{t}\,-\sum_{k=i+1}^{j} \frac{|t|^{m_kp}}{t^3}\,|\mu_k(0)|^p\,
			=\,\sum_{k=j+1}^{l} \frac{|t|^{m_kp}}{t^3}\,|\mu_k(0)|^p \\
			&\hspace{1cm}+ \left(\sum_{k=1}^{i} \gamma_k(t)\right) |t|^{q}\,+\sum_{k=i+1}^{l} |t|^{m_kp-2} \,\xi_k(t)+\Psi_{3}(t) +\mathcal{O}(|t|^{3q-3}),\quad t\in\Omega\setminus\{0\}.
		\end{align}
		Here again observe that the right hand side of \eqref{b8} is differentiable at $0$, which implies the differentiability of the left hand side of \eqref{b8}, and consequently there exists a non-empty set $J\subseteq \{i+1,\ldots,j\}$ such that $m_kp=2q$ for all $k\in J$ and $\alpha_{2}=\sum\limits_{k\in J}|\mu_k(0)|^p$. But $\alpha_{2}<0,$ and $\sum\limits_{k\in J}|\mu_k(0)|^p\geq 0$, which is a contradiction.
		
		\par Hence after analyzing all the above cases, we conclude that there is no isometric embedding of $S_q^m$ into $S_p^n$. This completes the proof.
	\end{proof}

	\begin{thm}\label{th3}
		Let $p\in(0,1)\setminus\{\frac{1}{k}:k\in\mathbb{N}\}$ and $2\leq m\leq n<\infty$. Then there is no isometric embedding of $S_1^m$ into $S_p^n$.
	\end{thm}
	\begin{proof}
		We will prove this result by contradiction. Suppose that there is an isometric embedding of $S_1^m$ into $S_p^n$. Here in the proof, we use the same notation that are used in the proof of Theorem \eqref{th1}. Therefore, we have the identity \eqref{a1} and consequently the identity \eqref{a2}, that is we have
		\begin{align}\label{c1}
			\left(1+|t|\right)^{p}=\sum_{k=1}^{n}|\lambda_k(t)|^p.
		\end{align}
		It is easy to observe that not all $\lambda_k(0)$ vanish. The non-differentiability of the left hand side of \eqref{c1} at $0$ forces that not all $\lambda_k(0)$ are non-zero. Therefore we may assume that there is some $l\in\{1,\ldots,n-1\}$ such that $\lambda_k(0)=0$ for $1\leq k\leq l$ and  $\lambda_k(0)\neq0$ for $l+1\leq k\leq n$. Hence we have the identity \eqref{a5}. That is 
		\begin{align}\label{c2}
			\nonumber&1+p\,|t|+\frac{p(p-1)}{2!}\,|t|^2+\mathcal{O}(|t|^3)=\sum_{k=1}^{l}|t|^{m_k p}|\,\mu_k(t)|^p+t\Psi_{1}(t)+1,\quad t\in\Omega,\\
			\implies& p\,|t|=\sum_{k=1}^{l}|t|^{m_k p}|\,\mu_k(t)|^p+t\Psi_{1}(t)-p(p-1)\,|t|^2+\mathcal{O}(|t|^3)\quad t\in\Omega.
		\end{align}
		Note that $m_kp\neq 1,\forall k\in\{1,\ldots,l\}$ as $p\in (0,1)\setminus\{\frac{1}{k}:k\in\mathbb{N}\}$. It is clear from the right hand side of \eqref{c2} that not all $m_kp>1$, otherwise the right hand side of \eqref{c2} is differentiable at $0$ but the left hand side is not differentiable at $0$. So there exists a non-empty subset $J\subseteq \{1,\ldots,l\}$ such that $m_kp<1$ for $k\in J$ and $m_kp>1$ for $k\in \{1,\ldots,l\}\setminus J$. Therefore, from \eqref{c2} we have
		\begin{align}\label{c3}
			\nonumber&p\,|t|-\sum_{k\in J}|t|^{m_k p}|\,\mu_k(0)|^p\\
			=&\sum_{k\in J}t|t|^{m_k p} \xi_k(t)+\sum_{k\in \{1,\ldots,l\}\setminus J}|t|^{m_k p} |\mu_k(t)|^p+t\Psi_{1}(t)-p(p-1)\,|t|^2+\mathcal{O}(|t|^3).
		\end{align}
		Note that the left hand side of \eqref{c3} is not differentiable at $0$ but the right hand side of \eqref{c3} is differentiable at $0$, which contradicts our assumption. This completes the proof.
	\end{proof}	
	
	\begin{thm}\label{th4}
		Let $p\in(0,1)\setminus\{\frac{1}{k}:k\in\mathbb{N}\}$ and $2\leq m\leq n<\infty$. Then there is no isometric embedding of $S_\infty^m$ into $S_p^n$.
	\end{thm}	
	
	\begin{proof}
		Suppose that there exists an isometric embedding $S_\infty^m$ into $S_p^n$. Then there exist $A$ and $B$ in $S_p^n$ such that $A, B$ have $(\mathbf{I}_{\infty,p})$. Therefore we have the following identity
		\begin{align}\label{d1}
			\max\{1,|t|^p\}=\|A+tB\|_p^p, \quad t\in\mathbb{R}.
		\end{align}
		By the Kato-Rellich Theorem \ref{kato}, there exist $n$ real-analytic functions $\lambda_1(\cdot), \lambda_2(\cdot),\ldots, \lambda_n(\cdot)$ in a small neighborhood $U$ of $1$ such that they are complete set of eigenvalues of $A+tB$ for $t\in U$. So for $t\in U $, from \eqref{d1} we have 
		\begin{align}\label{d2}
			\max\{1,|t|^p\}=\sum_{k=1}^{n}|\lambda_k(t)|^p.
		\end{align}
		If  $\lambda_k(1)\neq 0$ for all $1\leq k\leq n$, that is if $A+B$ is invertible then the right hand side of \eqref{d2} is a real-analytic  function in $U$ but the left hand side of \eqref{d2} is not analytic in $U$. So $A+B$ is a singular matrix. In \eqref{d2}, note that not all $\lambda_k(1)$ are zero as $1=\sum_{k=1}^{n}|\lambda_k(1)|^p$. So without loss of generality we may assume that there is a natural number $l\in\{1,2,\ldots,n\}$ such that $\lambda_k(1)=0$ for $1\leq k\leq l$ and  $\lambda_k(1)\neq 0$ for $l+1\leq k\leq n$. Note that if some $\lambda_i(t)$ is identically zero in $U$, we can neglect that one as it has no contribution in $\|A+tB\|_p$. Therefore we may also assume that all $\lambda_i(t)$ are non-zero analytic functions in $U$. Let $m_k$ be the multiplicity of zero at $1$ of the function $\lambda_k(t)$ for $1\leq k\leq l$. Let $\lambda_k(t)=(t-1)^{m_k}\mu_k(t), \mu_k(1)\neq 0,1\leq k\leq l$, for some analytic function $\mu_k(\cdot)$ in $U$. Then from \eqref{d2} we have,
		\begin{align}\label{d3}
			\max\{1,|t|^p\}=\sum_{k=1}^{l}|(t-1)|^{m_kp} |\mu_k(t)|^p+ \sum_{k=l+1}^{n} |\lambda_k(t)|^p,  \quad t\in U.
		\end{align}
		Since the left hand side of \eqref{d3} is not differentiable at $1$ so not all $m_kp$ in \eqref{d3} are greater than $1$. So we may assume that there is a  non-empty set $J\subseteq \{1,\ldots, l\}$ such that $m_kp<1$ for $k\in J$ and $m_kp>1$ for $k\in \{1,\ldots, l\}\setminus J$. Let $|\mu_k(t)|^p=|\mu_k(1)|^p+(t-1)\xi_k(t), k\in J$, where $\xi_k(\cdot)$ are some analytic function in $U$. Then from \eqref{d3}, we have
	\begin{align}\label{d4}
			\nonumber&\max\{1,|t|^p\}-\sum_{k\in J}|(t-1)|^{m_kp} |\mu_k(1)|^p\\
			=&\sum_{k\in J}|(t-1)|^{m_kp} (t-1)\,\xi_k(t)+ \sum_{k\in \{1,\ldots, l\}\setminus J}|(t-1)|^{m_kp} |\mu_k(t)|^p+\sum_{k=l+1}^{n} |\lambda_k(t)|^p,  \quad t\in U.
		\end{align}
		It is clear that the right hand side of \eqref{d4} is differentiable at $1$ but the left hand side is not differentiable at $1$, which leads to a contradiction. Therefore, there is no isometric embedding of $S_\infty^m$ into $S_p^n$. This completes the proof of the theorem.
	\end{proof}

	The following Lemma is needed to prove our main result in this section. 
	\begin{lma}\label{diflem}
		Let $0<q<\infty$, $0<p<1$, and $p\neq q$. Let $A$ and $B$ be two operators on $\ell_2^n$ such that $(A,B)$ has $(\mathbf{I}_{q,p})$. If $A<0$ or $A> 0$ then
		\begin{align*}
			\frac{d^2}{dt^2}\bigg|_{t=0}\norm{A+tB}_p^p< 0.
		\end{align*}
	\end{lma}

	\begin{proof}[Proof of Lemma~\ref{diflem}]
		Let $\{d_1,d_2,\ldots,d_n\}$ be the complete set of eigenvalues of $A$. Note that all $d_k, 1\leq k\leq n$ are non-zero and have the same sign. Since $A$ is invertible, there exists a real number $\epsilon>0$ such that $\sigma(A+tB)\cap [-\epsilon,\epsilon]=\emptyset$ for $t$ in a small neighbourhood of $0$. Now we define the compactly supported smooth function $f_p:\R\rightarrow \mathbb{C}$ such that $f_p(x)=|x|^p$ on $[-2,-\epsilon]\cup[\epsilon,2]$. Thus by applying Theorem~\ref{diffor} corresponding to the function $f_p$ we get 
		\begin{align}\label{e1}
			\nonumber\frac{d^2}{dt^2}\bigg|_{t=0}\norm{A+tB}_p^p=&\Tr\big(T^{A,A,A}_{f_p^{[2]}}(B,B)\big)\\
			=&\sum_{l=1}^{n}|b_{ll}|^2f_p^{[2]}(d_l,d_l,d_l)+\sum_{l=1}^{n-1}\sum_{k=l+1}^{n}~|b_{lk}|^2~\big(f_p^{[2]}(d_l,d_k,d_k)+f_p^{[2]}(d_k,d_l,d_l)\big).
		\end{align} 
		For the last identity in \eqref{e1}, we refer to \cite[Equation 4.3]{JFA22}.
		Note that $f_p^{[2]}(d_l,d_l,d_l)=p(p-1)|d_l|^{p-1}<0$  and for $d_l\neq d_k$, $\big(f_p^{[2]}(d_l,d_k,d_k)+f_p^{[2]}(d_k,d_l,d_l)\big)=p\cdot d_l^{p-2}\cdot\dfrac{1-\left(\frac{d_k}{d_l}\right)^{p-1}}{1-\left(\frac{d_k}{d_l}\right)}<0$.  Since $(A,B)$ has $(\mathbf{I}_{q,p})$, so by Corollary \ref{neq0}, we have $AB\neq 0$. Therefore along the same line of argument of the proof of \cite[Lemma 4.1]{JFA22}, from \eqref{e1}, we conclude the proof.
	\end{proof}
 Note that in the above Lemma \ref{diflem}, we prove that $\frac{d^2}{dt^2}\bigg|_{t=0}\norm{A+tB}_p^p< 0$ under some conditions on $A$ and $B.$ We want to point out that these type of results are often useful in various context. We refer to \cite{Zh20} where the author obtained an alternative proof of a famous theorem of \cite{BaCaLi94} (see also \cite{RiXu16} and \cite{Ko84}). To obtain this we use theory of multiple operator integrals and operator derivatives. We refer \cite{PoSkSu13}, \cite{PoSkSuOP13}, \cite{PoSu14} for more on this direction. 
 
\vspace{.2in} 
Now we are in a position to state and prove our main theorem in this section.
	\begin{thm}\label{th5}
		Let $0<p<1,\,2\leq q< \infty,\, 2\leq m\leq n<\infty$.  Then there is no isometric embedding $T:S_q^m\to S_p^n$ with $T(\text{diag}(1,0,\ldots,0))=A,$ and $T(\text{diag}(0,1,\ldots,0))=B$ such that 
		\begin{itemize}
			\item  $A,B\in M_n^{sa},$ 
			\item either $A\geq 0$ or $A\leq 0$.
		\end{itemize}
	\end{thm}
	\begin{proof}
		Here we use the same notations that are used in the proof of Theorem \ref{th1}. On the contrary, suppose there exists such an isometric embedding. Then we have the identity \eqref{a1}, and consequently \eqref{a2}. Altogether, we have
		\begin{align}\label{e2}
			\left(1+|t|^q\right)^{\frac{p}{q}}=\|A+tB\|_p^p=\sum_{k=1}^{n}|\lambda_k(t)|^p.
		\end{align}
		If $A$ is invertible, then the right hand side of \eqref{e2} is analytic in $\Omega$. Therefore, $\|A+tB\|_p^p$ has a power series representation in $\Omega$. Consequently from \eqref{e2} we have
		\begin{align}\label{e3}
			1+\frac{p}{q}|t|^q+\mathcal{O}(|t|^{2q})=1+t\frac{d}{dt}\bigg|_{t=0}\norm{A+tB}_p^p+\frac{t^2}{2!}\frac{d^2}{dt^2}\bigg|_{t=0}\norm{A+tB}_p^p+\mathcal{O}(|t|^{3q}),~~ t\in \Omega.
		\end{align}
		It is easy to observe that $\frac{d}{dt}\bigg|_{t=0}\norm{A+tB}_p^p=0$ as $1\leq \left(1+|t|^q\right)^{\frac{p}{q}}=\norm{A+tB}_p^p$   for all $t\in \Omega$. Therefore from \eqref{e3} we have
		\begin{align}\label{e4}
			\frac{p}{q}|t|^{q-2}+\mathcal{O}(|t|^{2q-2})=\frac{1}{2!}\frac{d^2}{dt^2}\bigg|_{t=0}\norm{A+tB}_p^p+\mathcal{O}(|t|^{3q-2}),~~ t\in \Omega\setminus\{0\}.
		\end{align}
		Taking limit $t\to 0$ on both sides of \eqref{e4}, we have
		\begin{align*}
			\frac{p}{q}~\lim_{t\to 0}|t|^{q-2}=\frac{1}{2!}\frac{d^2}{dt^2}\bigg|_{t=0}\norm{A+tB}_p^p< 0,
		\end{align*} which is impossible as $\frac{p}{q}>0$. We conclude that $A$ must be a singular matrix. Then, by the similar argument as given in the proof of Theorem~\ref{th2}, we can deduce the identity \eqref{a5} from the identity \eqref{a2}. Thus we have
		\begin{align}\label{e5}
			\frac{p}{q}|t|^{q-1}+\mathcal{O}(|t|^{2q-1})=\sum_{k=1}^{l}|t|^{m_k p-1}|\,\mu_k(t)|^p+\Psi_{1}(t).
		\end{align}
		Observe that $\Psi_{1}(0)=0$. Let $\Psi_{1}(t)=t\Psi_2(t)$, where $\Psi_{2}(\cdot)$ is some analytic function in $\Omega$. Then from \eqref{e5} we have,
		\begin{align}\label{e6}
			\frac{p}{q}|t|^{q-2}+\mathcal{O}(|t|^{2q-2})=&\sum_{k=1}^{l}|t|^{m_k p-2}|\,\mu_k(t)|^p+\Psi_{2}(t), \quad t\in \Omega;\\
			\implies \frac{p}{q}~\lim_{t\to 0}|t|^{q-2}= &\lim_{t\to 0}\,\sum_{k=1}^{l}|t|^{m_k p-2}|\,\mu_k(t)|^p +\Psi_{2}(0).
		\end{align}
		As $q\geq2$, from the above equation \eqref{e5}, it follows that  $m_kp\geq 2$ for all $k\in\{1,\ldots,l\}$. Let $J\subseteq\{1,\ldots,l\}$ be such that $m_kp=2$ for $k\in J$. Note that $J$ could be empty. Then by repeating the similar kind of arguments as given in \textbf{
			Subcase-II} of the proof of \cite[Theorem 4.2]{JFA22},  we have 
		\[\at{\frac{d^2}{dt^2}}{t=0}\norm{A+tB}_p^p=\sum_{k\in J}|\,\mu_k(0)|^p+\Psi_2(0)< 0.\quad (\text{ if } J=\emptyset, \text{ we assume }\sum_{k\in J}|\,\mu_k(0)|^p=0 )\]
		
		Therefore from \eqref{e6}, we conclude that
		
		\begin{align*}
			\frac{p}{q}~\lim_{t\to 0}|t|^{q-2}=\sum_{k\in J}|\,\mu_k(0)|^p+ \Psi_{2}(0)<0, \text{ but  } \frac{p}{q}>0,
		\end{align*}
		which leads to a contradiction. This completes the proof.
	\end{proof}
	\begin{crlre}
		Let $\mathcal M$ be a finite dimensional von Neumann algebra with normal semifinite faithful trace, and $\text{dim}\mathcal M\geq 2.$ Let $m\geq 2.$ Let $p\in(0,1)$. Then there is no isometric embedding of $S_q^m$ into $L_p(\mathcal M)$ for $q\in(0,2)\setminus\{1\}$, $S_1^m$ into $L_p(\mathcal M)$ for $p\in(0,1)\setminus\{\frac{1}{k}:k\in\mathbb N\}$ and  $S_\infty^n$ into $L_p(\mathcal M)$ for $p\in(0,1)\setminus\{\frac{1}{k}:k\in\mathbb N\}.$
	\end{crlre}
\begin{rmrk}
	The restrictions on the range of $p$ and $q$ in the above mentioned theorems is purely of technical nature. The main problem is that whenever we have nice analyticity of both the power-series in hand, we have no more tools to compare them as the eigenvalues of the analytic family $t\mapsto A+tB$ are in an abstract form if we use the Kato-Rellich theorem. Therefore, one needs some new ideas to solve the Question \ref{Xu} for the remaining cases.
\end{rmrk}
	\section*{Acknowledgements}  The first named author thankfully acknowledges the financial support provided by Mathematical Research Impact Centric Support (MATRICS) grant, File no: MTR/2019/000640, by the Science and Engineering Research Board (SERB), Department of Science and Technology (DST), Government of India. {\color{red}{The second named author acknowledges the financial support provided by NSF of China (No. 12071355) and the Fundamental Research Funds for the Central Universities (No. 2042022kf1185).}} The third named author thanks the Indian Institute of Technology Guwahati, Government of India, for the financial support. The fourth named author acknowledges {\color{red}{the financial support provided by}} DST-INSPIRE Faculty Fellowship No. DST/INSPIRE/04/2020/001132.

\end{document}